\documentclass[11pt]{article}

\usepackage{a4wide}

\usepackage{amsmath,amsbsy,amssymb,latexsym}
\usepackage{amssymb,amsthm}
\usepackage{url}

\newtheorem{theorem}{Theorem}[section]

\newtheorem{definition}[theorem]{Definition}
\newtheorem{remark}{Remark}

\begin{document}

\title{Calculus of variations with fractional derivatives\\
and fractional integrals\thanks{Accepted (July 6, 2009)
for publication in \emph{Applied Mathematics Letters}.}}

\author{Ricardo Almeida\\
\url{ricardo.almeida@ua.pt}
\and
Delfim F. M. Torres\\
\url{delfim@ua.pt}}

\date{Department of Mathematics,
University of Aveiro\\
3810-193 Aveiro, Portugal}

\maketitle


\begin{abstract}
We prove Euler-Lagrange fractional equations
and sufficient optimality conditions for problems of the
calculus of variations with functionals containing both fractional derivatives
and fractional integrals in the sense of Riemann-Liouville.

\bigskip

Riemann-Liouville fractional integral.

\smallskip

\end{abstract}


\section{Introduction}

In recent years numerous works have been dedicated
to the fractional calculus of variations.
Most of them deal with Riemann-Liouville fractional derivatives
(see \cite{Atanackovic,Baleanu,El-Nabulsi:Torres07,Frederico:Torres07,Frederico:Torres08}
and references therein), a few with Caputo or Riesz derivatives \cite{agrawal2,CD:Agrawal:2007,Almeida2,Baleanu:Agrawal}.
Depending on the type of functional being considered,
different fractional Euler-Lagrange type equations are obtained. We also mention \cite{muslih}, where a fractional Euler-Lagrange equation is obtained corresponding to a prescribed fractional space. Here we propose a new kind of functional with a
Lagrangian containing not only a Riemann-Liouville
fractional derivative (RLFD) but also
a Riemann-Liouville fractional integral (RLFI).
We prove necessary conditions of Euler-Lagrange type
for the fundamental fractional problem of the calculus of variations
and for the fractional isoperimetric problem.
Sufficient optimality conditions are also obtained
under appropriate convexity assumptions.


\section{Fractional Calculus}
\label{sec2}

Fractional calculus is an interdisciplinary area, with many applications in several fields, such as engineering \cite{ferreira,kulish,magin}, chemistry \cite{hilfer,metzler}, electrical and electromechanical systems \cite{debnath,Oustaloup}, viscoplasticity \cite{Diethelm}, physics \cite{hilfer2}, etc. \cite{TM}.

In this section we review the necessary definitions and facts from
fractional calculus. For more on the subject
we refer the reader to \cite{kilbas,Oldham,Podlubny,samko}.

Let $f$ be a function defined on the interval $[a,b]$.
Let $\alpha$ be a positive real and $n := [\alpha]+1$.

\begin{definition}
\label{def:frac:I:D}
The left RLFI is defined by
$${_aI_x^\alpha}f(x)=\frac{1}{\Gamma(\alpha)}\int_a^x (x-t)^{\alpha-1}f(t)dt,$$
and the right RLFI by
$${_xI_b^\alpha}f(x)=\frac{1}{\Gamma(\alpha)}\int_x^b(t-x)^{\alpha-1} f(t)dt.$$
The left RLFD is defined by
$${_aD_x^\alpha}f(x)=\frac{d^n}{dx^n} {_aI_x^{n-\alpha}}f(x)
=\frac{1}{\Gamma(n-\alpha)}\frac{d^n}{dx^n}\int_a^x
(x-t)^{n-\alpha-1}f(t)dt,$$
while the right RLFD is given by
$${_xD_b^\alpha}f(x)=(-1)^n\frac{d^n}{dx^n} {_xI_b^{n-\alpha}}f(x)
=\frac{(-1)^n}{\Gamma(n-\alpha)}\frac{d^n}{dx^n}\int_x^b
(t-x)^{n-\alpha-1} f(t)dt.$$
\end{definition}

The operators of Definition~\ref{def:frac:I:D}
are obviously linear. We now present the rules
of fractional integration by parts
for RLFI and RLFD. Let $p\geq1$, $q \geq 1$,
and $1/p+1/q\leq1+\alpha$. If $g\in L_p(a,b)$
and $f\in L_q(a,b)$, then
$$\int_{a}^{b}  g(x) \, {_aI_x^\alpha}f(x)dx
=\int_a^b f(x) \, {_x I_b^\alpha} g(x)dx \, ;$$
if $f$, $g$, and the fractional derivatives ${_aD_x^\alpha} g$
and ${_xD_b^\alpha} f$ are continuous on $[a,b]$,
then
$$ \int_{a}^{b}g(x) \, {_aD_x^\alpha} f(x)dx
=\int_a^b f(x) \, {_xD_b^\alpha} g(x)dx\, ,$$
$0<\alpha<1$.

\begin{remark}
\label{Bound2}
The left RLFD of $f$ is infinite at $x=a$ if $f(a) \not=0$ (cf. \cite{ross}). Similarly, the right RLFD is infinite
if  $f(b)\not=0$. Thus, assuming that $f$ possesses  continuous
left and right RLFD on $[a,b]$, then $f(a)=f(b)=0$
must be satisfied.
\end{remark}


\section{The Euler-Lagrange equation}
\label{sec3}

Let us consider the following problem:
\begin{equation}
\label{Funct1}
\mathcal{J}(y)=\int_a^b
L(x,{_a I_x^{1-\alpha}} y(x),\,{_aD_x^\beta}y(x)) \, dx
\longrightarrow \min.
\end{equation}
We assume that
$L(\cdot,\cdot,\cdot) \in C^1([a,b]\times\mathbb{R}^2; \mathbb{R})$,
$x\to \partial_2 L(x,{_a I_x^{1-\alpha}} y(x),\,{_aD_x^\beta}y(x))$
has continuous right RLFI of order $1-\alpha$ and
$x\to \partial_3 L(x,{_a I_x^{1-\alpha}} y(x),\,{_aD_x^\beta}y(x))$
has continuous right RLFD of order $\beta$, where
$\alpha$ and $\beta$ are real numbers in the interval $(0,1)$.

\begin{remark}
We are assuming that the admissible functions $y$ are such that
${_a I_x^{1-\alpha}} y(x)$ and ${_aD_x^\beta}y(x)$ exist
on the closed interval $[a,b]$. We also note that
as $\alpha$ and $\beta$ goes to $1$ our fractional functional
$\mathcal{J}$ tends to the classical functional
$\int_a^b L(x,y(x),y'(x)) \, dx$
of the calculus of variations.
\end{remark}

\begin{remark}
We consider functionals $\mathcal{J}$ containing
the left RLFI and the left RLFD only.
This comprise the important cases in applications.
The results of the paper are easily generalized for functionals
containing also the right RLFI and/or right RLFD.
\end{remark}

\begin{theorem}[The fractional Euler-Lagrange equation]
\label{Theo E-L1} Let $y(\cdot)$ be a local minimizer of problem
(\ref{Funct1}).
Then, $y(\cdot)$ satisfies the fractional Euler-Lagrange equation
\begin{equation}
\label{E-L1}
{_x I_b^{1-\alpha}} \partial_2L(x,{_a I_x^{1-\alpha}} y(x),\,{_aD_x^\beta}y(x))
+ {_xD_b^\beta}\partial_3L(x,{_a I_x^{1-\alpha}} y(x),\,{_aD_x^\beta}y(x)) = 0
\end{equation}
for all $x\in[a,b]$.
\end{theorem}

\begin{remark}
Condition \eqref{E-L1} is only necessary for an extremum.
The question of sufficient conditions for an extremum
is considered in Section~\ref{sec4}.
\end{remark}

\begin{proof}
Since $y$ is an extremizer of $\mathcal{J}$,
by a well known result of the calculus of variations the first variation of
$\mathcal{J}(\cdot)$ is zero at $y$, \textrm{i.e.},
\begin{equation}
\label{eq1}
0
= \displaystyle \delta \mathcal{J}(\eta,y)
=\displaystyle \int_a^b({_a I_x^{1-\alpha}} \eta \, \partial_2L
+ {_aD_x^\beta}\eta \, \partial_3L  )\, dx \, .
\end{equation}
Integrating by parts,
\begin{equation}
\label{IP1}
\displaystyle\int_a^b {_a I_x^{1-\alpha}} \eta \,
\partial_2L \, dx = \int_a^b \eta \, {_x I_b^{1-\alpha}} \partial_2L \, dx
\end{equation}
and
\begin{equation}
\label{IP2}
\displaystyle
\int_a^b {_aD_x^\beta}\eta \, \partial_3L \, dx
=\int_a^b \eta \, {_xD_b^\beta}\partial_3L\, dx.
\end{equation}
Substituting (\ref{IP1}) and (\ref{IP2}) into equation (\ref{eq1}),
we find that
$\displaystyle \int_a^b({_x I_b^{1-\alpha}}
\partial_2L+ {_xD_b^\beta}\partial_3L  )\eta\, dx=0$
for each $\eta$. Since $\eta$ is an arbitrary function,
by the fundamental lemma of the calculus of variations we deduce that
${_x I_b^{1-\alpha}} \partial_2L+ {_xD_b^\beta}\partial_3L =0$.
\end{proof}

\begin{remark}
As $\alpha$ and $\beta$ goes to $1$,
the fractional Euler-Lagrange equation (\ref{E-L1}) becomes
the classical Euler-Lagrange equation $\partial_2L-d/dx\partial_3L=0$.
\end{remark}

A curve that is a solution of the fractional differential equation
(\ref{E-L1}) will be called an \emph{extremal} of $\mathcal{J}$.
Extremals play also an important role in the solution
of the fractional isoperimetric problem
(see Section~\ref{sec6}). We note that equation (\ref{E-L1}) contains
right RLFI and right RLFD, which are not present in the formulation
of problem (\ref{Funct1}).


\section{Some generalizations}
\label{sec5}

We now give some generalizations of Theorem \ref{Theo E-L1}.

\subsection{Extension to variational problems of non-commensurate order}

We now consider problems of the calculus of variations with Riemann-Liouville derivatives and integrals of non-commensurate order, \textrm{i.e.},
we consider functionals containing RLFI and RLFD of different fractional orders. Let
\begin{equation}
\label{Funct2}
\mathcal{J}(y)=\int_a^b L(x,{_a I_x^{1-\alpha_1}} y(x),\ldots,{_a I_x^{1-\alpha_n}} y(x),{_aD_x^{\beta_1}}y(x),\ldots,{_aD_x^{\beta_m}}y(x)) \, dx,
\end{equation}
where $n$ and $m$ are two positive integers and $\alpha_i,\beta_j\in(0,1)$, $i=1,\ldots,n$ and $j=1,\ldots,m$.
Following the proof of Theorem \ref{Theo E-L1}, we deduce the following result.

\begin{theorem}
If $y(\cdot)$ is a local minimizer of (\ref{Funct2}), then $y(\cdot)$ satisfies the Euler-Lagrange equation
$$\sum_{i=1}^n{_x I_b^{1-\alpha_i}} \partial_{i+1}L + \sum_{j=1}^m{_xD_b^{\beta_j}}\partial_{j+n+1}L=0$$
for all $x\in[a,b]$.
\end{theorem}

\subsection{Extension to several dependent variables}

We now study the case of multiple unknown functions $y_1,\ldots,y_n$.

\begin{theorem}
Let $\mathcal{J}$ be the functional given by the expression
$$\mathcal{J}(y_1,\ldots,y_n)=\int_a^b L(x,{_a I_x^{1-\alpha}}y_1(x),\ldots,{_a I_x^{1-\alpha}} y_n(x),{_aD_x^{\beta}}y_1(x),\ldots,{_aD_x^{\beta}}y_n(x)) \, dx.$$
If $y_1(\cdot),\ldots,y_n(\cdot)$ is a local minimizer of $\mathcal{J}$, then it satisfies for all $x\in[a,b]$ the following system of $n$ fractional differential equations:
$${_xI_b^{1-\alpha}} \partial_{k+1}L+{_xD_b^{\beta}}\partial_{n+k+1}L=0, \quad k=1,\ldots,n.$$
\end{theorem}

\begin{proof}
Denote by $y$ and $\eta$ the vectors $(y_1,\ldots,y_n)$ and $(\eta_1,\ldots,\eta_n)$, respectively.
For a parameter $\epsilon$, we consider a new function
\begin{equation}\label{eq2}
J(\epsilon)=\mathcal{J}(y+\epsilon \eta) \, .
\end{equation}
Since $y_1(\cdot),\ldots,y_n(\cdot)$ is an extremizer of $\mathcal{J}$, $J'(0)=0$. Differentiating equation (\ref{eq2}) with respect to $\epsilon$, at $\epsilon=0$, we obtain
$$\displaystyle \int_a^b\left[{_a I_x^{1-\alpha}} \eta_1 \,\partial_2L+\cdots + {_a I_x^{1-\alpha}} \eta_n \,\partial_{n+1}L+ {_aD_x^\beta}\eta_1 \,
\partial_{n+2}L +\cdots+ {_aD_x^\beta}\eta_n \,\partial_{2n+1}L \right]\, dx=0.$$
Integrating by parts leads to
$$\displaystyle \int_a^b\left[{_xI_b^{1-\alpha}} \partial_2L + {_xD_b^\beta}\partial_{n+2}L\right]\,\eta_1+\cdots + \left[{_xI_b^{1-\alpha}} \partial_{n+1}L
+ {_xD_b^\beta}\partial_{2n+1}L \right]\,\eta_n \, dx=0.$$
Considerer a variation $\eta = (\eta_1,0,\ldots,0)$, $\eta_1$ arbitrary; then by the fundamental lemma of the calculus of variations we obtain
${_xI_b^{1-\alpha}} \partial_2L+ {_xD_b^\beta}\partial_{n+2}L=0$.
Selecting appropriate variations $\eta$, one deduce the remaining formulas.
\end{proof}


\section{The fractional isoperimetric problem}
\label{sec6}

We consider now the problem of minimizing the functional $\mathcal{J}$ given by (\ref{Funct1}) subject to an integral constraint
$\mathcal{I}(y)=\int_a^b g(x,{_a I_x^{1-\alpha}} y(x),\,{_aD_x^\beta}y(x)) \, dx=l$,
where $l$ is a prescribed value. This problem was solved in \cite{Almeida2} for functionals containing Caputo fractional derivatives and RLFI.
Using similar techniques as the ones discussed in \cite{Almeida2}, one proves the following:

\begin{theorem}\label{thm:nc:ip}
Consider the problem of minimizing the functional $\mathcal{J}$ as in (\ref{Funct1}) on the set of functions $y$ satisfying condition
 $\mathcal{I}(y)=l$. Let $y$ be a local minimum for the problem. Then, there exist two constants $\lambda_0$ and $\lambda$, not both zero,
such that $y$ satisfies the Euler-Lagrange equation
${_x I_b^{1-\alpha}} \partial_2K+ {_xD_b^\beta}\partial_3K=0$ for all $x \in [a,b]$,
where $K=\lambda_0 L+\lambda g$.
\end{theorem}

\begin{remark}
If $y$ is not an extremal for $\mathcal{I}$, then one can choose $\lambda_0 = 1$ in Theorem~\ref{thm:nc:ip}:
there exists a constant $\lambda$ such that $y$ satisfies ${_x I_b^{1-\alpha}} \partial_2F+ {_xD_b^\beta}\partial_3F=0$ for all  $x \in [a,b]$, where $F=L+\lambda g$.
\end{remark}


\section{Sufficient conditions}
\label{sec4}

In this section we prove sufficient conditions that ensure the existence of minimums. Similarly to what happens in the classical calculus of variations,
some conditions of convexity are in order.

\begin{definition}
Given a function $L$, we say that $L(\underline x,u,v)$ is convex in $S\subseteq\mathbb R^3$ if $\partial_2 L$ and $\partial_3 L$
exist and are continuous and verify the following condition:
$$L(x,u+u_1,v+v_1)-L(x,u,v)\geq \partial_2 L(x,u,v)u_1+\partial_3 L(x,u,v)v_1$$
for all $(x,u,v),(x,u+u_1,v+v_1)\in S$.
\end{definition}

Similarly, we define convexity for $L(\underline x, \underline u,v)$.

\begin{theorem}
Let $L(\underline x,u,v)$ be a convex function in $[a,b]\times\mathbb R^2$ and let $y_0$ be a curve satisfying the fractional Euler-Lagrange equation (\ref{E-L1}).
Then, $y_0$ minimizes (\ref{Funct1}).
\end{theorem}

\begin{proof}
The following holds:
$$\begin{array}{ll}
\mathcal{J}(y_0+\eta)-\mathcal{J}(y_0)&=
\displaystyle \int_a^b \left[L(x,{_a I_x^{1-\alpha}} y_0(x)
+{_a I_x^{1-\alpha}} \eta(x),{_aD_x^\beta}y_0(x)+{_aD_x^\beta}\eta(x))\right.\\
&\left. \qquad- L(x,{_a I_x^{1-\alpha}} y_0(x),\,{_aD_x^\beta}y_0(x)) \right] dx \\
& \geq \displaystyle \int_a^b \left[ \partial_2 L(x,{_a I_x^{1-\alpha}} y_0(x),\,
{_aD_x^\beta}y_0(x))\, {_aI_x^{1-\alpha}}\eta \right.\\
&\left.\qquad+\partial_3 L(x,{_a I_x^{1-\alpha}} y_0(x),\,
{_aD_x^\beta}y_0(x))\,{_aD_x^\beta}\eta \right] \, dx\\
&= \displaystyle\int_a^b \left[{_xI_b^{1-\alpha}}\partial_2 L
+{_xD_b^\beta}\partial_3 L\right]_{ (x,{_a I_x^{1-\alpha}} y_0(x),\,
{_aD_x^\beta}y_0(x))}\eta\,dx=0.\\
\end{array}$$
Thus, $\mathcal{J}(y_0+\eta)\geq \mathcal{J}(y_0)$.
\end{proof}

We now present a sufficient condition for convex Lagrangians
on the third variable only. First we recall the notion of exact field.

\begin{definition} Let $D\subseteq \mathbb R^2$ and let $\Phi:D \to \mathbb R$ be a function of class $C^1$. We say that $\Phi$ is an exact field for $L$ covering $D$
if there exists a function $S \in C^1(D,\mathbb R)$ such that
$$\begin{array}{rl}
\partial_1S(x,y)&=L(x,y,\Phi(x,y))-\partial_3L(x,y,\Phi(x,y))\Phi(x,y)\, ,\\
\partial_2S(x,y)&=\partial_3L(x,y,\Phi(x,y))\, .\\
\end{array}$$
\end{definition}

\begin{remark}
This definition is motivated by the classical Euler-Lagrange equation. Indeed, every solution $y_0 \in C^2[a,b]$ of the differential equation
$y'=\Phi(x,y(x))$ satisfies the (classical) Euler-Lagrange equation $\partial_2L-\frac{d}{dx}\partial_3L=0$.
\end{remark}

\begin{theorem}
Let $L(\underline x, \underline u, v)$ be a convex function
in $[a,b]\times \mathbb R^2$, $\Phi$ an exact field for $L$ covering
$[a,b]\times \mathbb R \subseteq D$,
and $y_0$ a solution of the fractional equation
\begin{equation}
\label{FracEquation}
{_aD_x^\alpha}y(x)=\Phi(x,{_aI_x^{1-\alpha}}y(x)).
\end{equation}
Then, $y_0$ is a minimizer for
$\mathcal{J}(y)=\int_a^b L(x,{_a I_x^{1-\alpha}}
y(x),\,{_aD_x^\alpha}y(x)) \, dx$
subject to the constraint
\begin{equation}
\label{Constraints}
\left\{ y:[a,b]\to\mathbb R \, | \, {_aI_a^{1-\alpha}}y(a)
={_aI_a^{1-\alpha}}y_0(a), \,
{_aI_b^{1-\alpha}}y(b)={_aI_b^{1-\alpha}}y_0(b)\right\}.
\end{equation}
\end{theorem}

\begin{proof}
Let $E(x,y,z,w)=L(x,y,w)-L(x,y,z)-\partial_3L(x,y,z)(w-z)$.
First observe that
$$\begin{array}{ll}
\displaystyle\frac{d}{dx}S(x,{_a I_x^{1-\alpha}}y(x))
&=\displaystyle \partial_1S(x,{_a I_x^{1-\alpha}}y(x))
+\partial_2S(x,{_a I_x^{1-\alpha}}y(x))\displaystyle\frac{d}{dx}{_a I_x^{1-\alpha}}y(x)\\
&= \displaystyle\partial_1S(x,{_a I_x^{1-\alpha}}y(x))
+\partial_2S(x,{_a I_x^{1-\alpha}}y(x)){_a D_x^{\alpha}}y(x)\, .\\
\end{array}$$
Since $E \geq 0$, it follows that
$$\begin{array}{ll}
\mathcal{J}(y) &= \displaystyle \int_a^b\left[ E(x,{_a I_x^{1-\alpha}} y,\,
\Phi(x, {_a I_x^{1-\alpha}} y), \,{_aD_x^\alpha}y)
+ L(x,{_a I_x^{1-\alpha}} y, \Phi(x, {_a I_x^{1-\alpha}}y))\right.\\
&\qquad \left. + \displaystyle \partial_3L(x,{_a I_x^{1-\alpha}} y,
\Phi(x, {_a I_x^{1-\alpha}}y))({_aD_x^\alpha}y
-\Phi(x, {_a I_x^{1-\alpha}}y))\right]\, dx\\
&\geq \displaystyle  \int_a^b\left[ L(x,{_a I_x^{1-\alpha}} y,
\Phi(x, {_a I_x^{1-\alpha}}y))\right.\\
&\left.\qquad + \displaystyle  \partial_3L(x,{_a I_x^{1-\alpha}} y,
\Phi(x, {_a I_x^{1-\alpha}}y))({_aD_x^\alpha}y-\Phi(x, {_a I_x^{1-\alpha}}y))\right]\, dx\\
&= \displaystyle \int_a^b \left[ \partial_1S(x, {_a I_x^{1-\alpha}}y)
+ \partial_2S(x, {_a I_x^{1-\alpha}}y) {_aD_x^\alpha}y  \right] \, dx\\
&=\displaystyle  \int_a^b \frac{d}{dx}S(x, {_a I_x^{1-\alpha}}y) \, dx\\
&=S(b, {_a I_b^{1-\alpha}}y(b))- S(a, {_a I_a^{1-\alpha}}y(a)) .\\
\end{array}$$
Because $y_0$ is a solution of (\ref{FracEquation}),
$E(x,{_a I_x^{1-\alpha}} y_0,\,\Phi(x, {_a I_x^{1-\alpha}} y_0), \,{_aD_x^\alpha}y_0)=0$.
With similar calculations as before, one has
$\mathcal{J}(y_0)=S(b,{_a I_b^{1-\alpha}}y_0(b))- S(a, {_a I_a^{1-\alpha}}y_0(a))$.
We just proved that $\mathcal{J}(y_0)\leq \mathcal{J}(y)$ when subject
to the constraint (\ref{Constraints}).
\end{proof}


\section{Conclusions}

In this note we consider a new class of fractional functionals of the calculus of variations that depend not only on fractional derivatives but also on fractional integrals.
We exhibit necessary and sufficient conditions of optimality
for the fundamental problem of the calculus of variations
and for problems subject to integral constrains (isoperimetric problems). As future work it would be interesting to address the question of existence of solutions, and to study direct methods
to minimize the proposed type of functionals.


\section*{Acknowledgments}

The authors are grateful to the support
of the \emph{Control Theory Group}
from the \emph{Centre for Research on Optimization and Control}
given by the \emph{Portuguese Foundation for Science and Technology} (FCT), cofinanced through the
\emph{European Community Fund} FEDER/POCI 2010.



\end{document}